\documentclass[12pt]{amsart}

\usepackage{times,epsf,amssymb,amsmath,hyperref, yhmath, wrapfig, rlepsf, wasysym,color}

\begin{document}

\newtheorem{thm}{Theorem}[section]
\newtheorem{lem}[thm]{Lemma}
\newtheorem{cor}[thm]{Corollary}
\newtheorem{conj}[thm]{Conjecture}

\theoremstyle{definition}
\newtheorem{defn}[thm]{\bf{Definition}}

\theoremstyle{remark}
\newtheorem{rmk}[thm]{Remark}

\theoremstyle{question}
\newtheorem{question}[thm]{Question}

\def\square{\hfill${\vcenter{\vbox{\hrule height.4pt \hbox{\vrule width.4pt height7pt \kern7pt \vrule width.4pt} \hrule height.4pt}}}$}

\newcommand{\SI}{\partial_\infty ({\Bbb H}^2\times {\Bbb R})}
\newcommand{\Si}{S^1_{\infty}\times {\Bbb R}}
\newcommand{\si}{S^1_{\infty}}
\newcommand{\CS}{S^1_\infty\times\{\pm\infty\}}
\newcommand{\PI}{\partial_{\infty}}
\newcommand{\caps}{{\Bbb H}^2\times\{\pm\infty\}}

\newcommand{\BH}{\Bbb H}
\newcommand{\BHH}{{\Bbb H}^2\times {\Bbb R}}
\newcommand{\BR}{\Bbb R}
\newcommand{\BC}{\Bbb C}
\newcommand{\BZ}{\Bbb Z}

\newcommand{\wh}{\widehat}
\newcommand{\wt}{\widetilde}
\newcommand{\wc}{\check}

\newcommand{\A}{\mathcal{A}}
\newcommand{\bb}{\mathcal{B}}
\newcommand{\Oo}{\mathcal{O}}
\newcommand{\V}{\mathcal{V}}
\newcommand{\U}{\mathcal{U}}
\newcommand{\X}{\mathcal{X}}
\newcommand{\Y}{\mathcal{Y}}
\newcommand{\W}{\mathcal{W}}
\newcommand{\Z}{\mathcal{Z}}
\newcommand{\C}{\mathcal{C}}
\newcommand{\K}{\mathcal{K}}
\newcommand{\D}{\mathcal{D}}
\newcommand{\p}{\mathcal{P}}
\newcommand{\Q}{\mathcal{Q}}
\newcommand{\R}{\mathcal{R}}
\newcommand{\h}{\mathcal{H}}
\newcommand{\pp}{\mathcal{P}}
\newcommand{\T}{\mathcal{T}}
\newcommand{\N}{\mathcal{N}}
\newcommand{\s}{\mathcal{S}}

\newcommand{\e}{\epsilon}
\newcommand{\B}{\mathbf{B}}

\newenvironment{pf}{{\it Proof:}\quad}{\square \vskip 12pt}

\title[Minimal Surfaces in $\BHH$]{Minimal Surfaces in $\BHH$: Nonfillable Curves}
\author{Baris Coskunuzer}
\address{UT Dallas, Dept. Math. Sciences, Richardson, TX 75080}
\email{baris.coskunuzer@utdallas.edu}
\thanks{The author is partially supported by Simons Collaboration Grant, and Royal Society Newton Mobility Grant.}

\maketitle

\begin{abstract}

In this paper, we study the asymptotic Plateau problem in $\BHH$. We give the first examples of non-fillable finite curves with no thin tail in $\Si$. Furthermore, we study the fillability question for infinite curves in $\SI$.
\end{abstract}

\section{Introduction}

Asymptotic Plateau Problem in $\BHH$ asks the existence of a minimal surface $\Sigma$ in $\BHH$ for a given curve $\Gamma$ in $\SI$ with $\PI\Sigma=\Gamma$. We will call a finite collection of disjoint Jordan curves $\Gamma$ in $\SI$ {\em fillable}, if $\Gamma$ bounds a complete, embedded, {\em minimal surface}  $S$ in $\BHH$ with $\PI S=\Gamma$. We will call $\Gamma$ {\em strongly fillable} if $\Gamma$ bounds a complete, embedded, {\em area minimizing surface} $\Sigma$ in $\BHH$ with $\PI \Sigma=\Gamma$. In our previous papers \cite{Co1, Co2}, we finished the  classification of strongly fillable curves in $\SI$. In this paper, we will study the fillable curves in $\SI$.

In the past decade, existence, uniqueness, and regularity of minimal and CMC surfaces in $\BHH$ have been studied extensively, and many important results have been  obtained by the leading researchers of the field, e.g. \cite{CHR, CR, CMT, Da, FMMR, HRS, KM, MMR, MoR, MRR, NR, PR, RT, ST, ST2}.

Let $\Gamma$ be a finite collection of Jordan curves in $\SI$. Define
$\Gamma^\pm=\Gamma\cap (\overline{\BH^2}\times\{\pm\infty\})$. We call $\Gamma$ {\em an infinite curve} if either $\Gamma^+$ or $\Gamma^-$ is nonempty. We will call $\Gamma$ {\em a finite curve} otherwise. In particular, $\Gamma$ is a finite curve if $\Gamma$ belongs to the open cylinder at infinity, i.e. $\Si$. 

In \cite{ST}, Sa Earp and Toubiana gave an obstruction for $\Gamma$ in $\SI$ to be fillable, i.e. having a thin tail (See Lemma \ref{thinlem}). This is the only obstruction known for a finite curve to be fillable so far. Hence, the following question is quite natural:

\begin{question} Let $\Gamma$ be a finite collection of pairwise disjoint finite Jordan curves in $\Si$. Then, does $\Gamma$ bound an embedded minimal surface in $\BHH$ if $\Gamma$ has no thin tail?
\end{question}

In this paper, we will answer this question negatively by constructing the first examples of finite non-fillable curves with no thin tail. 

\begin{thm} There exists finite non-fillable curves with no thin tail in $\Si$.	
\end{thm}

On the other hand, we will discuss fillability question for infinite curves in $\SI$, and give many fillable, and non-fillable examples. 


The organization of the paper is as follows. In the next section, we give some definitions, and overview the related results on the fillability and strong fillability questions. In Section 3, we will construct the new families of non-fillable finite curves with no thin tail. In Section 4, we will study the fillable and non-fillable infinite curves in $\SI$. Finally, in Section 5, we will discuss further questions, and give some concluding remarks.

\subsection{Acknowledgements}

I would like to thank Rafe Mazzeo, Giuseppe Tinaglia, and Francisco Martin for very valuable remarks.

\section{Preliminaries}

In this section, we will give the basic definitions, and a brief overview of the past results which will be used in the paper. 

Throughout the paper, we use the product compactification of $\BHH$. In particular, $\overline{\BHH}=\overline{\BH^2}\times\overline{\BR}=\BHH\cup \SI$ where $\SI$ consists of three components, i.e. the infinite open cylinder $\Si$ and the closed caps at infinity $\overline{\BH^2}\times\{+\infty\}$, $\overline{\BH^2}\times\{-\infty\}$. Hence, $\overline{\BHH}$ is a solid cylinder under this compactification.

Let $\Sigma$ be an open, complete surface in $\BHH$, and $\PI\Sigma$ represent the asymptotic boundary of $\Sigma$ in $\SI$. Then, if $\overline{\Sigma}$ is the closure of $\Sigma$ in $\overline{\BHH}$, then $\PI \Sigma= \overline{\Sigma}\cap \SI$.

\begin{defn} A surface is {\em minimal} if the mean curvature $H$ vanishes everywhere. A compact surface with boundary $\Sigma$ is called {\em area minimizing surface} if $\Sigma$ has the smallest area among the surfaces with the same boundary. A noncompact surface is called {\em area minimizing surface} if any compact subsurface is an area minimizing surface. 
\end{defn}

In this paper, we will study the Jordan curves in $\SI$ which bounds a complete, embedded, minimal surfaces in $\BHH$.
Throughout the paper, when we say {\em a curve in $\SI$} we mean a finite collection of pairwise disjoint Jordan curves in $\SI$.

\begin{defn} [Fillable Curves] Let $\Gamma$ be a curve in $\SI$. We will call $\Gamma$ {\em fillable} if $\Gamma$ bounds a complete, embedded, minimal surface $S$ in $\BHH$, i.e. $\PI S=\Gamma$. We will call $\Gamma$ {\em strongly fillable} if $\Gamma$ bounds a complete, embedded, area minimizing surface $\Sigma$ in $\BHH$, i.e. $\PI \Sigma=\Gamma$. We call such $S$ or $\Sigma$ as {\em filling surface} for $\Gamma$. 
\end{defn}

Notice that a strongly fillable curve is fillable since any area minimizing surface is minimal. \textit{The asymptotic Plateau problem for $\BHH$} is the following classification problems:

\vspace{.2cm}

{\em Which $\Gamma$ in $\SI$ is fillable or strongly fillable?}

\vspace{.2cm}

Throughout the paper, we will use the following notation for a curve $\Gamma$ in $\SI$. Decompose $\Gamma=\Gamma^+\cup\Gamma^-\cup\wt{\Gamma}$ such that
$\Gamma^\pm=\Gamma\cap (\overline{\BH^2}\times\{\pm\infty\})$ and $\wt{\Gamma}=\Gamma\cap(\Si)$. In particular, $\Gamma^\pm$ is a collection of closed arcs and points in the closed caps at infinity, where $\wt{\Gamma}$ is a collection of open arcs and closed curves in the infinite open cylinder. With this notation, we will call a curve $\Gamma$ {\em finite} if $\Gamma^+=\Gamma^-=\emptyset$. In other words, $\Gamma$ is a finite curve if $\Gamma$ belongs to the open cylinder at infinity $\Si$. We will call $\Gamma$ {\em infinite} otherwise. In this part, we will overview the basic results for finite curves. 

\subsection {Thin Tails} \

One of the most interesting properties of the asymptotic Plateau problem in $\BHH$ is the existence of non-fillable curves. While any curve $\Lambda$ in $S^2_\infty(\BH^3)$ is strongly fillable in $\BH^3$ \cite{An}, Sa Earp and Toubiana showed that there are some non-fillable $\Gamma$ in $\SI$ \cite{ST}.

\begin{defn} \label{taildef} [Thin tail] Let $\Gamma$ be a Jordan curve in $\SI$, and let $\tau$ be an arc in $\Gamma$. Assume that there is a vertical straight line $L_0$ in $\Si$ such that

\begin{itemize}

\item $\tau \cap L_0 \neq \emptyset$ and $\partial \tau\cap L_0 = \emptyset$,

\item $\tau$ stays in one side of $L_0$,

\item $\tau\subset \si \times (c,c+\pi)$ for some $c\in \BR$.

\end{itemize}

Then, we call $\tau$ {\em a thin tail} in $\Gamma$.

\end{defn}

\begin{lem} \cite{ST} \label{thinlem} Let $\Gamma$ be a curve in $\SI$. If $\Gamma$ contains a thin tail, then there is no properly immersed minimal surface $\Sigma$ in $\BHH$ with $\PI \Sigma=\Gamma$.
\end{lem}

\subsection{Tall Curves} \

The above result shows that the curves with thin tail cannot be fillable. Hence, to bypass this obstruction, we introduced the following notion in \cite{Co1}.

\begin{defn} \label{tallrecdef} [Tall Rectangles] \cite{Co1} Consider $\Si$ with the coordinates $(\theta, t)$ where $\theta\in [0,2\pi)$ and $t\in \BR$. We will call the rectangle $R=[\theta_1,\theta_2]\times[t_1,t_2]\subset \Si$  {\em tall rectangle} if $t_2-t_1> \pi$.
\end{defn}

\begin{rmk} \label{tallrecrmk}
The boundaries of tall rectangles bounds a special minimal surfaces which can be described explicitly \cite{ST}. In particular, let $R_h$ be a rectangle of height $h=t_2-t_1>\pi$. Let $\tau_h=\partial R_h$. Then, by \cite{Co1}, $\tau_h$ bounds a unique minimal surface, which is a minimal plane $\pp_h$ in $\BHH$ with $\PI \pp_h=\tau_h$. Furthermore, by \cite{ST}, when $h\nearrow \infty$, $\pp_h$ gets closer to the vertical geodesic plane $\gamma\times\BR$ where $\gamma$ is the geodesic in $\BH^2$ connecting $\theta_1$ and $\theta_2$ in $\si$. When $h\searrow\pi$, then $\pp_h$ completely escapes to infinity. See section \cite[Section 2.2]{Co1} for further details.
\end{rmk}


\begin{defn} \label{talldef} [Tall Curves]
We call a curve $\Gamma$ in $\SI$  {\em tall curve} if the open region $\Si-\Gamma$ can be covered by tall rectangles $int(R_i)$, i.e. $\Si-\Gamma=\bigcup_i int(R_i)$. 
\end{defn}

Notice that if $\Gamma$ is a tall curve, then it cannot contain a thin tail. Note also that this definition naturally generalizes to infinite curves \cite{Co2}. In \cite{Co1}, we also defined the {\em height of a curve}, $h(\Gamma)$, as the length of the smallest component in vertical line segments in $\Si-\Gamma$. Hence, $\Gamma$ is tall if and only if $h(\Gamma)>\pi$. We will call a curve $\Gamma$ {\em short} if $h(\Gamma)<\pi$. 

Note that throughout the paper, we will not consider the curves with height $h(\Gamma)=\pi$ as it is not hard to construct fillable, or non-fillable curves with $h(\Gamma)=\pi$, e.g. rectangle of height $\pi$ (nonfillable), and a butterfly curve with neck height $\pi$ (fillable).

By using this notion, we gave a fairly complete classification of \textit{strongly fillable finite curves} in \cite{Co1} as follows.

\begin{lem} \cite{Co1} \label{finiteAPP} Let $\Gamma$ be a finite curve in $\Si$ with $h(\Gamma)\neq \pi$. Then, $\Gamma$ is strongly fillable if and only if $\Gamma$ is a tall curve.
\end{lem}

A reasonable question at this point is whether {\em tall condition} is equivalent to {\em no thin tail condition} for a curve $\Gamma$ in $\SI$. The answer to this question is "No". When $\Gamma$ is not connected, or null-homotopic in $\Si$, there are trivial examples for this question as described below. When $\Gamma$ is null-homotopic Jordan curve in $\Si$, the examples are more delicate: Butterfly curves (explained below).

\subsection{Fillable Short Curves} \label{fillshortsec} \

By Lemma \ref{finiteAPP}, any strongly fillable curve must be tall. However, this does not mean any fillable curve must be tall, too. By using a simple observation, it is easy to produce many examples of fillable short curves. In particular, union of tall curves may not be tall. Notice that fillability question does not enforce the minimal surface to be connected. So, if $\Gamma=\gamma_1\cup\gamma_2 ..\cup\gamma_n$ is a union of pairwise disjoint strongly fillable curves in $\si$, then $\Gamma$ is also fillable. One can see this as follows: Each $\gamma_i$ bounds an area minimizing surface $\Sigma_i$, i.e. $\PI\Sigma_i=\gamma_i$. As each $\Sigma_i$ is area minimizing, and $\{\gamma_i\}$ is pairwise disjoint, by \cite{Co1}, $\{\Sigma_i\}$ is pairwise disjoint, too. Hence, $\Sigma=\bigcup \Sigma_i$ is an embedded minimal surface with $\PI\Sigma=\Sigma_i$. 
	
Now, if we choose $\{\gamma_i\}$ vertically close to each other, then $\Gamma$ is no longer a tall curve. By using this idea, it is straightforward to obtain many fillable short curves, e.g. $\Gamma=\si\times\{c_1,c_2\}$ where $|c_1-c_2|<\pi$. Note that $\Sigma$ is only a minimal surface, not an area minimizing surface since any area minimizing surface must be bounded by a tall curve by Lemma \ref{finiteAPP}. So, such a $\Gamma$ is fillable, but not strongly fillable.

\vspace{.2cm}

\begin{wrapfigure}{r}{1.2in}
	
	\relabelbox  {\epsfxsize=1in
		
		\centerline{\epsfbox{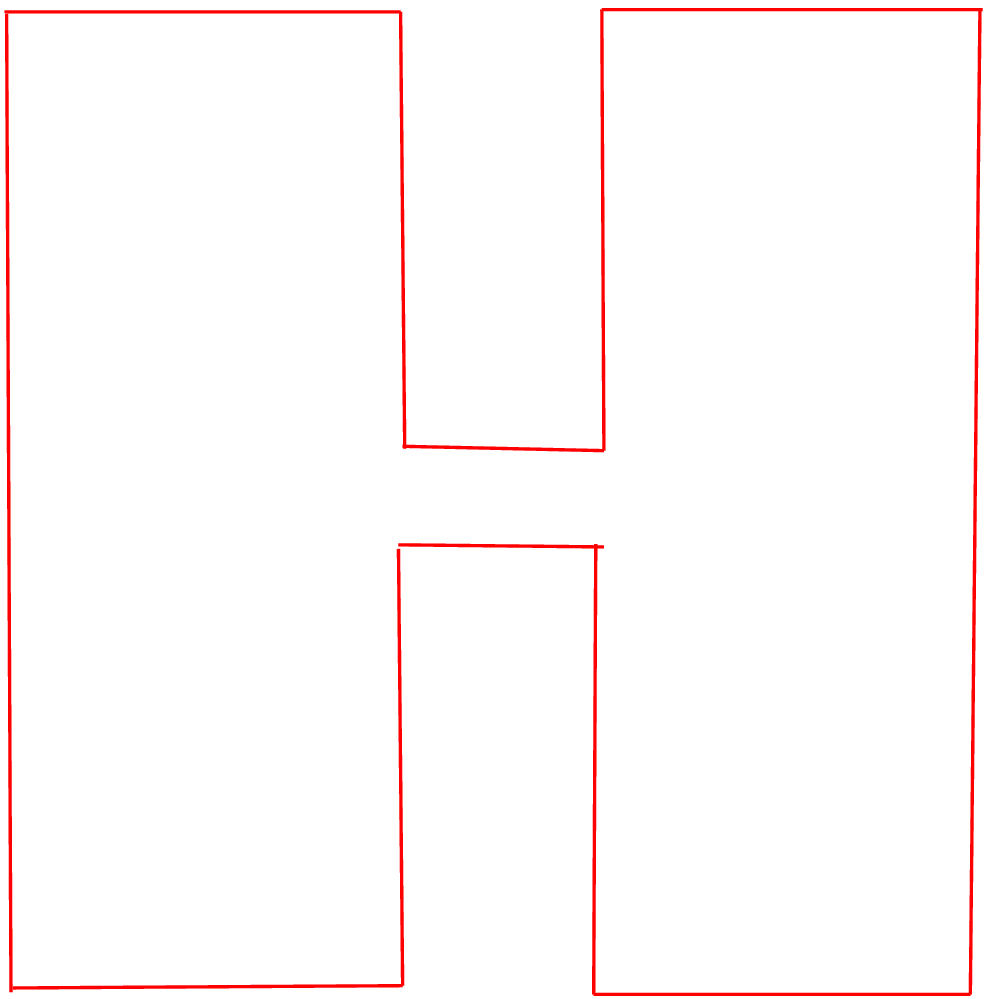}}}

	\endrelabelbox
	
\end{wrapfigure}

\noindent {\bf Butterfly Curves:} The previous short fillable examples came from arranging different components of a curve vertically close to each other. On the other hand, there is also a special family of \textit{connected} short fillable curves. In \cite{KM,Co1}, the authors constructed these null-homotopic curves in $\Si$ in the shape of a butterfly as in the figure right. In these curves, the neck height is $<\pi$ while the height of the rectangles is $>>\pi$. These curves are called {\em Butterfly Curves}. These  are important examples to see that "no thin tail" condition is not equivalent to being tall for a connected curve.

Another important family of fillable curves is the essential curves in $\Si$ which are graphs over $\si$. We will call such curves {\em graphical curves}.

\begin{lem} \label{graphlem} [Graphical Minimal Disks] Let $\gamma$ be a Jordan curve in $\Si$ which is graph of a smooth function $u:\si \to \BR$. Then, $\gamma$ bounds a unique minimal surface $\Sigma$ in $\BHH$. Furthermore, $\Sigma$ is the graph of a function $\wh{u}:\BH^2\to \BR$.
\end{lem}

\subsection{Vertical Catenoids} \label{catenoidsec} \

We will describe one more special family of minimal surfaces in $\BHH$ which will be one of the key components of our technique in the next section: {\em Rotationally symmetric vertical catenoids}. These are  minimal catenoids which are asymptotic to two horizontal round circles in $\Si$ \cite{NSST}. In particular, for $d\in(0,\pi)$, $\C_d$ is the minimal catenoid with $\PI \C_d=\si\times\{c,c+d\}$. As being rotationally symmetric, they are generated by rotating an infinite catenary arc about a vertical line in $\BHH$.

The behavior of the vertical catenoid $\C_d$ is completely determined by its height $d$. Let $\eta_d$ be the necksize of the vertical catenoid $\C_d$. Then, when $d\nearrow \pi$, then $\eta_d\to\infty$, and when $d\searrow 0$, then $\eta_d\to 0$. In other words, when $d\nearrow\pi$, the catenoid $\C_d$ gets closer and closer to the asymptotic sphere $\Si$, and escapes to infinity. On the other hand, when $d\searrow 0$, $\C_d$ limits to the double cover of a horizontal geodesic plane $\BH^2\times\{c\}$.	See \cite[Section 7.1]{Co1} for explicit definition of $\C_d$ and further details.

\section{Non-fillable Finite Curves}

As mentioned in the introduction, so far, the only known finite non-fillable curves in $\Si$ are the ones with thin tail (Definition \ref{taildef}). In this section, we will construct the first examples of non-fillable curves with no thin tail. 

Our curves are a generalization of the butterfly curves defined earlier. In \cite{Co1, KM}, the authors showed the existence of fillable short curves by using these butterfly curves. Here, we will generalize their definition, and get non-fillable curves with no thin tail. 

\begin{figure}[h]
	
	\relabelbox  {\epsfxsize=3.5in
		
		\centerline{\epsfbox{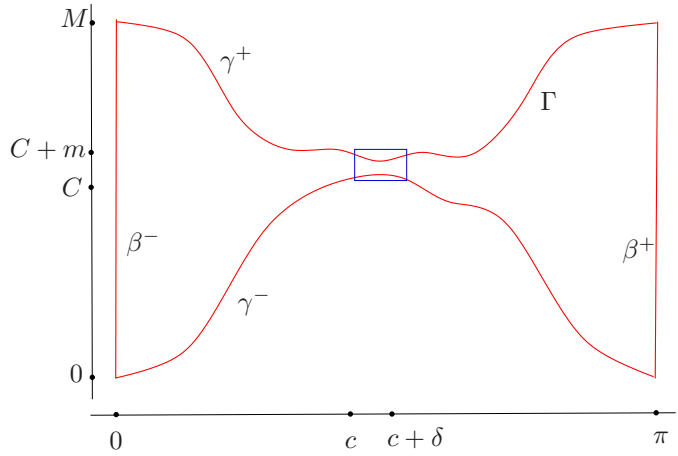}}}
	
	\relabel{1}{\footnotesize $M$}
	\relabel{2}{\footnotesize $C+m$}
	\relabel{3}{\footnotesize $C$}
	\relabel{4}{\footnotesize $0$}
	\relabel{5}{\footnotesize $0$}
	\relabel{6}{\footnotesize $c$}
	\relabel{7}{\footnotesize $c+\delta$}
	\relabel{8}{\footnotesize $\pi$}
	\relabel{9}{\footnotesize $\gamma^+$}
	\relabel{10}{\footnotesize $\gamma^-$}
	\relabel{11}{\footnotesize $\beta^-$}
	\relabel{12}{\footnotesize $\beta^+$}
	\relabel{13}{\footnotesize $\Gamma$}

	\endrelabelbox
	
	\caption{\label{genbutfig} \small $\Gamma$ is a $(m,\delta,M)$-butterfly curve in $\Si$. $\Gamma$ is in a rectangle of height $M$, and contains a thin segment in a rectangle of dimensions $\delta \times m$.}
	
\end{figure}

\begin{defn} \label{genbutdefn} [Generalized Butterfly Curves] Let $0<m<\pi$,  $\delta>0$, and $M>\pi$. Let $u^\pm:[0,\pi]\to [0,M]$ be smooth functions such that 

\begin{enumerate}
	\item $u^+(\theta)>u^-(\theta)$ for $\theta\in[0,\pi]$
	\item $u^+(0)=u^+(\pi)=M$ and $u^-(0)=u^-(\pi)=0$
	\item For some $c\in (0,\pi-\delta)$, and $C\in(0,M-m)$, $u^\pm([c,c+\delta])\subset [C,C+m]$. 	
\end{enumerate}	

Let $\alpha^\pm$ be the arcs in $\Si$ which are the graphs of $u^\pm$.  Let $\beta^-=\{0\}\times[0,M]$ and $\beta^+=\{\pi\}\times[0,M]$ be the vertical line segments in $\Si$. Define $\Gamma=\alpha^+\cup\alpha^-\cup\beta^+\cup\beta^-$. We call $\Gamma$ a  \textit{$(m,\delta,M)$-butterfly curve} (See Figure \ref{genbutfig}).
\end{defn}

Notice that as $M>\pi$, any such generalized butterfly curve has no thin tail. Note also that even though the above definition looks very restrictive, the proof of the theorem below show that it can be generalized to a rich family of curves holding basic conditions. See Remark \ref{genbutrmk}.


To prove our nonexistence theorem, we need the following lemma. Let $R$ be a tall rectangle as in Definition \ref{taildef}, i.e.  $R=[\theta_1,\theta_2]\times [C_1,C_2]$ be a rectangular region in $\Si$ with $C_2-C_1>\pi$. Let $\tau=\partial R$. Let $\pp$ be the unique minimal plane in $\BHH$ with $\partial \pp=\tau$ (Remark \ref{tallrecrmk}). Let $\Delta$ be the component of $\BHH-\pp$ with $\PI\Delta=R$.

\begin{lem} \label{traplem} Let $R, \tau, \pp, \Delta$ be as above. Let $\Gamma$ be a Jordan curve in $R\subset \Si$. If $\Sigma$ is a minimal surface in $\BHH$ with $\PI \Sigma=\tau$, then $\Sigma\subset \Delta$.	
\end{lem}

\begin{pf} Let $\wh{\theta}_1=\theta_1-\e$ and $\wh{\theta}_2=\theta_2+\e$ for some small $\e>0$. Consider the geodesic $\eta$ in $\BH^2$ connecting $\wh{\theta}_1$ and $\wh{\theta}_2$. Consider the geodesic plane $\eta\times \BR$ in $\BHH$. Let $U$ be the component of $\BHH-\eta\times\BR$ containing $\Delta$. As $U$ is convex in $\BHH$, $\Sigma\subset U$.
	
Now, we will foliate $U-\Delta$ by minimal planes spanning tall rectangles. Let $\rho:[0,\epsilon)\to[0,\infty)$ be a monotone increasing, onto function. Let $R_t=[\theta_1-t,\theta_1+t]\times[C_1-\rho(t), C_2+\rho(t)]$ for $t\in [0,\e)$. Let $\tau_t=\partial R_t$. Let $\pp_t$ be the unique area minimizing surface in $\BHH$ with $\PI\pp_t=\tau_t$. Then, by \cite[Section 2.2]{Co1}, $\pp_t\cap\pp_s=\emptyset$ for $t\neq s$, and $\{\pp_t\mid t\in [0,\e)\}$ is a foliation of $U-\Delta$.

Now, by above $\Sigma\subset U$. If $\Sigma$ is not completely in $\Delta$, then let $t_0=\sup\{t\mid \pp_t\cap \Sigma\neq \emptyset\}$. Then, $\pp_{t_0}$ and $\Sigma$ has a tangential intersection where one lies in one side of the other. This contradicts to the maximum principle.
\end{pf}

\begin{thm} For any $m_0<\pi$ and $\delta_0>0$, there exists an $M_0=M(m_0,\delta_0)>\pi$ such that for any $m<m_0$, $\delta>\delta_0$ and $\pi<M<M_0$, any $(m,\delta,M)$-butterfly curve is non-fillable.	
\end{thm}

\begin{pf} Let $\Gamma_M$ be such a $(m,\delta,M)$-butterfly curve. For a given $m_0<\pi$ and $\delta_0>0$, we will determine $M_0=M(m_0,\delta_0)>\pi$ later. Let $\Sigma_M$ be an embedded minimal surface in $\BHH$ with $\PI\Sigma_M=\Gamma_M$. Our strategy is first to build a cage $\Omega_M$ for $\Sigma_M$. Then, we will run over  $\Sigma_M$ trapped in the cage $\Omega_M$ with a vertical catenoid, and get a contradiction by maximum principle. 

\vspace{.2cm}

\noindent {\bf Step 1:} Construction of the "thin" cage $\Omega_M$ for $\Sigma_M$, i.e. $\Sigma_M\subset \Omega_M$.

\vspace{.2cm}

\noindent {\em Proof of Step 1:} Consider the rectangle $R_M=[0,\pi]\times[0,M]$ in $\Si$. Then, $\Gamma_M\subset R_M$. Let $\tau_M=\partial R_M$, and $\pp_M$ be the minimal plane in $\BHH$ with $\PI\pp_M=\tau_M$. Let $\Delta_M$ be the component of $\BHH-\pp_M$ with $\PI\Delta_M=R_M$. Then, as $\Gamma_M\subset R_M$, $\Sigma_M\subset \Delta_M$ by Lemma \ref{traplem}.

Now, consider the functions $u_M^\pm:[0,\pi]\to[0,M]$ in the construction of the butterfly curve as in Definition \ref{genbutdefn}. Extend $u_M^\pm$ to $\wh{u}_M^\pm:\si\to \BR$ trivially such that $\wh{u}_M^+(\theta)=M$ and $\wh{u}_M^-(\theta)=0$ for any $\theta \in [\pi,2\pi]$. Let $\wh{\gamma}_M^\pm$ be Jordan curves in $\Si$ which are graphs of functions $u_M^\pm$. By \cite{CR}, $\wh{\gamma}^\pm_M$ bounds a unique minimal surface $\Q_M^\pm$ which is a graph over $\BH^2$. Let $V_M$ be the region between $\Q_M^+$ and $\Q_M^-$ in $\BHH$. Then, $\Sigma_M$ must be in $V_M$. In order to see this, foliate the region above $\Q_M^+$ by $\Q^+_{Mt}=\Q_M+t$ for $t>0$, i.e. vertical translation of $\Q_M$ by $t$. Then, if $\Sigma_M$ goes over $\Q_M^+$, it must intersect some $\Q^+_{Mt}$. Let $t_0=\inf\{t\mid\Sigma_M\cap\Q^+_{Mt}\neq \emptyset\}$. Minimal surfaces $\Q^+_{Mt_0}$ and $\Sigma_M$ have tangential intersection by lying in one side. This contradicts to the maximum principle. Similarly, $\Sigma_M$ cannot go below $\Q_M^-$. Hence, $\Sigma_M\subset V_M$.

Let $\Omega_M=\Delta_M\cap V_M$. As $\Sigma_M\subset \Delta_M$ and $\Sigma_M\subset V_M$, $\Sigma_M\subset \Omega_M$.  We will call $\Omega_M$ a {\em cage for $\Sigma_M$} (See Figure \ref{cage}). Step 1 follows.

\begin{figure}[h]
	
	\relabelbox  {\epsfxsize=2.5in
		
		\centerline{\epsfbox{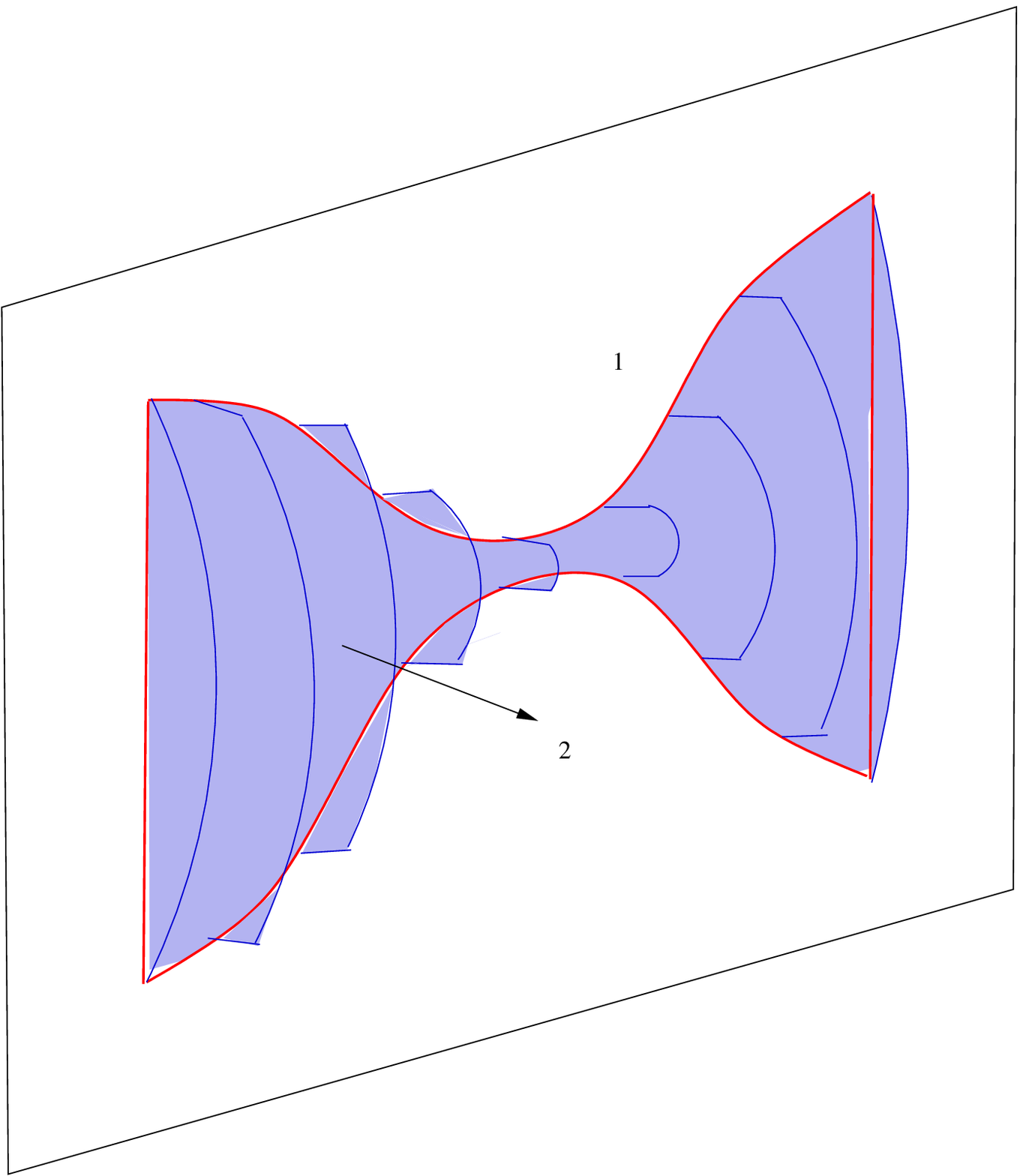}}}
	
	\relabel{1}{\footnotesize $\Gamma_M$}
	\relabel{2}{\footnotesize $\Omega_M$}

	\endrelabelbox
	
	\caption{\label{cage} \small For a given $\Gamma_M$ in $\Si$, we build a cage $\Omega_M$ to contain any minimal surface $\Sigma_M$ in $\BHH$ with $\PI\Sigma_M=\Gamma_M$.}
	
\end{figure}

\vspace{.2cm}

\noindent {\bf Step 2:} Killing $\Sigma_M$ with a vertical catenoid.

\vspace{.2cm}

\noindent {\em Proof of Step 2:} Recall that by Definition \ref{genbutdefn},  for some $c\in (0,\pi-\delta_0)$, and $C\in(0,M-m_0)$,  $u_M^\pm([c,c+\delta_0])\subset [C,C+m_0]$. Fix sufficiently small $\e>0$ such that $4\e<\pi-m_0$. Consider the rotationally symmetric vertical catenoid $\C_h$ of height $h=m_0+4\e$ such that $\PI\C_h=\si\times\{C-2\e,C+m_0+2\e\}$ (Section \ref{catenoidsec}). Let $k^+=C+m_0+\e$ and $k^-=C-\e$. Consider the compact catenoid $\wh{\C}_h=\C_h\cap \BH^2\times [k^-,k^+]$. Our aim is to translate $\wh{\C}_h$ with hyperbolic isometries towards $\Sigma_M$ so that first point of touch in the interior gives a contradiction with maximum principle. To implement this, we need to make sure that the boundary of the translation of $\wh{\C}_h$ always stays away from $\Omega_M$, and hence away from $\Sigma_M$. We will achieve this by choosing $M$ sufficiently close to $\pi$ so that the cage $\Omega_M$ is very "thin", i.e. $\Omega_M$ is "very close" to infinity. 

Let $\mu$ be the geodesic in $\BH^2$ asymptotic to the points $c+\frac{\delta}{2}$ and $\frac{3\pi}{2}$. Let $\varphi_s$ be the hyperbolic isometry fixing $\mu$ pushing from $\frac{3\pi}{2}$ to $c+\frac{\delta}{2}$ with translation length $s>0$. For $s<0$, let $\varphi_s=\varphi^{-1}_{|s|}$. Let $\wt{\varphi}_s$ be the isometry of $\BHH$ with $\wt{\varphi}_s(x,t)=(\varphi_s(x),t)$. Consider the horizontal tranlation of compact catenoid $\wh{\C}_h$ by $\wt{\varphi}_s$, i.e. $\wh{\C}_h^s=\wt{\varphi}_s(\wh{\C}_h)$. Let $Y_h^+$ be the disk in $\BH^2\times\{k^+\}$ bounded by the round circle $\partial^+ \wh{\C}_h$, and let $r_h$ be the radius of $Y_h$. Let $Y_h^{s+}= \wt{\varphi}_s (\Y_h^+) $. Define $\bb^+_{r_h}=\bigcup_{s\in\BR} Y_h^{s+}$. By construction, $\partial \bb^+_{r_h}$ will be a pair (Right and Left) of $r_h$-equidistant curves $\tau^+_R$ and $\tau^+_L$ to $\mu$ in $\BH^2\times \{k^+\}$ (See Figure \ref{banana}). Similarly, define $\bb^-_{r_h}$. $\partial \bb^-_{r_h}$ will be another pair of equidistant curves $\tau^-_R$ and $\tau^-_L$ in $\BH^2\times \{k^-\}$. Hence, $\bb^\pm_{r_h}$ are banana shapes in $\BH^2\times\{k^\pm\}$  with $\partial \bb^\pm_{r_h}=\tau^\pm_R\cup\tau^\pm_L$. The bottom tip of the bananas $\bb^\pm_{r_h}$ will be $(c+\frac{\delta}{2}, k^\pm)\in \Si$ by construction (See Figure \ref{banana}). 

\begin{figure}[h]
	
	\relabelbox  {\epsfxsize=2.5in
		
		\centerline{\epsfbox{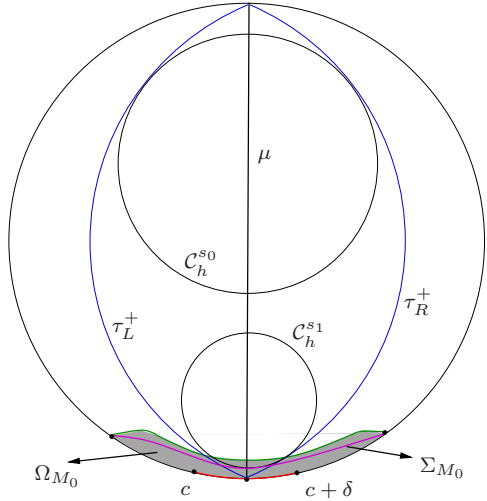}}}
	
	\relabel{1}{\tiny $\tau^+_R$}
	\relabel{2}{\tiny $\tau^+_L$}
	\relabel{3}{\tiny $\C_h^{s_0}$}
	\relabel{4}{\tiny $\C_h^{s_1}$}
	\relabel{5}{\tiny $c$}
	\relabel{6}{\tiny $c +\delta$}
	\relabel{7}{\tiny $\mu$}
	\relabel{8}{\tiny $\Omega_{M_0}$}
	\relabel{9}{\tiny $\Sigma_{M_0}$}

	\endrelabelbox
	
	\caption{\label{banana} \footnotesize The banana region $\bb^+_{r_h}$ is in $\BH^2\times\{C+m_0+\e\}$, and bounded by the blue equidistant curves. The shaded region  $\Omega_{M_0}$ is in $\BH^2\times \{C+m_0\}$ (level below), and near the red segment $[c,c+\delta]\times\{C+m_0\}$ in $\Si$. For $M_0$ sufficiently close to $\pi$, $\bb^+_{r_h}$ will be away from the cage $\Omega_{M_0}$. }
	
\end{figure}

We want to make sure that $\bb^\pm_{r_h}$ is disjoint from the cage $\Omega_M$. Notice that as $M\searrow \pi$ then $\pp_M \to \Si$, i.e. $\pp_M$ escapes to infinity (See Remark \ref{tallrecrmk}). In other words, the "distance" of $\pp_M$ to $\Si$ goes to $0$. Hence, as $M\searrow \pi$, the horizontal thickness of the cage $\Omega_M$ goes to $0$. Recall that $\delta_0$ is fixed, and $\Q_M^+$ is graphical near $[c,c+\delta_0]\times \{C+m_0\}$ by \cite{CR}. Therefore for sufficiently small $M_0>\pi$, the banana $\bb^+_{r_h}$ inside the higher level $\BH^2\times \{C+m_0+\e\}$ will be away from the cage $\Omega_{M_0}$ (See Figure \ref{banana}). Similarly, it is true for $\bb^-_{r_h}$. In particular, we have $\bb^\pm_{r_h}\cap\Omega_{M_0}=\emptyset$.

Now, consider $\wh{\C}_h^s=\wt{\varphi}_s(\wh{\C}_h)$. By construction, for sufficiently small $s_0<<0$, $\wh{\C}_h^{s_0}\cap \Omega_{M_0}=\emptyset$. As $\Sigma_{M_0}\subset \Omega_{M_0}$, $\wh{\C}_h^{s_0}\cap \Sigma_{M_0}=\emptyset$. Then, assuming $\Sigma_{M_0}$ is not empty, let $s_1=\inf\{s\mid \wh{\C}_h^s\cap\Sigma_{M_0}\neq \emptyset\}$.  Notice that we made sure that $\partial\wh{\C}_h^s\cap \Sigma_{M_0}=\emptyset$  for any $s\in \BR$, as $\bb^\pm_{r_h}\cap\Omega_{M_0}=\emptyset$. Then, we have $\wh{\C}_h^{s_1}$ and $\Sigma_{M_0}$ has a tangential interior intersection with one lies in one side of the other. As both are minimal surfaces, this gives a contradiction. 

Notice that if we take $m<m_0$, and $\delta>\delta_0$, we can still use the same cage $\Omega_{M_0}$ and the same vertical catenoid $\C_h$ of height $h=m_0+2\e$. The proof follows.	
\end{pf}

\begin{rmk} \label{genbutrmk} [Modification of Generalized Butterfly Curves] Notice that, in the proof of theorem above, the restriction for the domain of $u^\pm$ to be $[0,\pi]$ is indeed not essential. This is because for any proper subinterval $[\theta_1,\theta_2]$ in $\si$, by using hyperbolic isometries, we can make an isometric change in the setup for $\Gamma$ to satisfy the above $[0,\pi]$ condition. This condition is given in this way so that the constant $M_0=M(m_0,\delta_0)$ does not depend any other parameter. Clearly, any vertical translation does not change the setup either. Hence, as long as $\Gamma$ contains a thin strip in a rectangle of dimension $\delta\times m$, and it is contained in a rectangle of height $M(m,\delta)>\pi$, then $\Gamma$ will be a non-fillable curve with no thin tail.
\end{rmk}

\section{Fillable and Non-fillable Infinite Curves} \label{nonfillsec}

In this part, we will see that fillability question is quite different for infinite curves. In particular, for a given infinite curve $\Gamma$, we will show that $\Gamma^\pm$ may not be very useful to detect whether $\Gamma$ is fillable by constructing examples and non-examples. We will start with the notation and the background.

\subsection{Infinite Curves:} \label{infcurvesec} \

As we defined before if $\Gamma^\pm=\Gamma\cap (\overline{\BH^2}\times\{\pm\infty\}$) is nonempty, we call $\Gamma$ in $\SI$ an \textit{infinite curve}. In this paper, we will adapt the notation in \cite{Co2} for infinite curves. For further details, and examples, see \cite{Co2}.

We will decompose $\Gamma^\pm$ into two parts: $\Gamma^\pm=\overline{\Gamma}^\pm_g \cup \Gamma^\pm_c$. Let $\Gamma_g^\pm=\Gamma\cap \BH^2\times\{\pm\infty\}$ corresponds to the interior of $\Gamma^\pm$. In particular, if nonempty, $\Gamma_g^\pm$ will be a collection of arcs in $\BH^2\times\{\pm\infty\}$. 

Let $\overline{\Gamma_g^\pm}$ be the closure of $\Gamma_g^\pm$ in $\overline{\BH^2}\times\{\pm\infty\}$. Let $\Gamma_c^\pm=\Gamma^\pm-\overline{\Gamma_g^\pm}$ be the remaining boundary points of $\Gamma^\pm$ in $S^1_\infty\times\{\pm\infty\}$. Here, the subscript $c$ correspond to the term "corner". We will call $\Gamma$ {\em nonoverlapping at the corner} if $\Gamma_c^\pm$ does not contain any interval in the corner circles $\CS$. In other words, $\Gamma$ does not overlap with the corner circles of $\SI$ in any interval. 

We will call an infinite curve $\Gamma$ in $\SI$ {\em tame} if $\Gamma^\pm$ has finitely many components. Otherwise, we will call $\Gamma$ {\em a wild curve.} Note that throughout the paper, all curves are assumed to be tame unless otherwise stated.  

\vspace{.2cm}

\noindent {\em Infinite Rectangles:} Let $\gamma$ be a complete geodesic in $\BH^2$ with $\PI\gamma=\{p,q\}$. Let $\alpha$ be one of the two arcs in $S^1_\infty(\BH^2)$ with endpoints $p$ and $q$. Fix $t_0\in\BR$. Let $l^+_p=\{p\}\times[t_0,\infty]$ and $l^-_p=\{p\}\times[-\infty,t_0]$ be the vertical line segments in $\Si$. Let $\gamma^\pm=\gamma\times\{\pm\infty\}$ be the geodesic in $\caps$. Let $\alpha_0=\alpha \times\{t_0\}$. Then define $\R^+=\gamma^+\cup l^+_p\cup l^+_q\cup \alpha_0$ is an infinite rectangle (See Figure \ref{trappedfig}-left). Similarly, $\R^-=\gamma^-\cup l^-_p\cup l^-_q\cup \alpha_0$ is also an infinite rectangle. Note that any infinite rectangle $\R$ in $\SI$ bounds a unique area minimizing surface $\T$ in $\BHH$ with $\PI\T=\R$ \cite{Co2}.

\subsection{Fillability of Infinite Curves} \

Now, we will look at the properties of infinite fillable curves. A good starting point for our study would be the classification result for strongly fillable curves from \cite{Co2}.

\begin{lem} \cite{Co2} \label{stronglem} Let $\Gamma$ be a tame infinite curve in $\SI$. Then, $\Gamma$ is strongly fillable if and only if all of the following conditions satisfied:
	\begin{enumerate}
		\item $int(\Gamma^\pm)$ is a collection of geodesics (possibly empty).
		\item $\Gamma$ is tall.	
		\item $\Gamma$ is fat at infinity.
		\item $\Gamma$ is nonoverlapping at the corner.
	\end{enumerate}
\end{lem}

Notice that any strongly fillable curve is fillable. In particular, any curve holding these 4 conditions are automatically fillable. On the other hand, if $\Gamma$ is strongly fillable infinite curve, it must hold these 4 conditions. Hence, by considering above classification result, a good question would be the following:

\begin{question} Which of the properties above is also necessary for fillable infinite curves?
\end{question}

In other words, we want to understand if $\Gamma$ is fillable, which of these 4 properties must hold.

\vspace{.2cm}

\noindent {\bf (1) $\Gamma^\pm_g$ being geodesics:} By the following lemma, $\Gamma^\pm_g$ being geodesics is also a necessary condition for fillable curves. 

\begin{lem} \cite{KM} \label{geodlem} If $\Gamma$ is an infinite fillable curve in $\SI$, then $\Gamma_g^\pm$ must be a collection of geodesics in $\caps$.
\end{lem}

\vspace{.2cm}

\noindent {\bf (2) $\Gamma$ being tall:} Notice that being a tall curve is only related with the finite part $\wt{\Gamma}$ of the infinite curve $\Gamma$. Recall that being tall is not a necessary condition for being fillable, even though it is necessary condition for being strongly fillable. Similarly, we can produce many fillable short examples by using the short examples in Section \ref{fillshortsec}. Take any collection of these short fillable curves, and add an infinite fillable curve (like an infinite rectangle). Then, we would trivially have short infinite curves. One can also use two infinite rectangles whose finite horizontal segments are vertically $h<\pi$ close to get such curves.

One can ask if $\Gamma$ is infinite \textit{connected} curve, is being tall still necessary? Then, recall the finite butterfly curves in Section \ref{fillshortsec}, which are examples for connected fillable short curves. It is straightforward to generalize this construction to infinite curves. One can use two disjoint infinite geodesic planes, and connect them with a horizontal thin strip (height $h<\pi$) just like finite butterfly curves. Then, by using the barrier construction in \cite[Theorem 5.1]{Co1}, one can get connected,  short, fillable, infinite curve.

In particular, while being tall is a necessary condition for being strongly fillable, it is not necessary for being fillable.

\vspace{.2cm}

\noindent {\bf (3) $\Gamma$ being fat at infinity:} We refer to \cite{Co2} for definitions of \textit{fat at infinity}, and \textit{skinny at infinity} for a curve in $\SI$. Being {\em fat at infinity} is not a necessary condition for being fillable by the following lemma. 

\begin{lem} \cite[Theorem 4.3]{Co2} \label{fillem} Let $\gamma_1^+\cup..\cup\gamma_n^+$ and $\gamma_1^-\cup..\cup\gamma_m^-$ be given collections of pairwise disjoint geodesics in $\BH^2\times\{+\infty\}$ and $\BH^2\times\{-\infty\}$ respectively. Then, there exists a fillable curve $\Gamma$ with $\Gamma^+=\gamma_1^+\cup..\cup\gamma_n^+$ and $\Gamma^-=\gamma_1^-\cup..\cup\gamma_m^-$.
\end{lem}

This lemma shows that for fillability question, being fat or skinny at infinity is irrelevant. This is because for any collection of disjoint geodesics $\Lambda^\pm$ in $\BH^2\times\{\pm\infty\}$, it is possible to find a fillable curve $\Gamma$ with $\Gamma^\pm=\Lambda^\pm$. This means while being fat at infinity is a necessary condition for being strongly fillable, it is not necessary for being fillable (See Figure \ref{trappedfig}-left).

\vspace{.2cm}

\noindent {\bf (4) $\Gamma$ being nonoverlapping at the corner:} Recall that $\Gamma$ being nonoverlapping at the corner means that $\Gamma$ contains an interval in the corner circles $\si\times\{\pm\infty\}$. We show that this is indeed a necessary condition for being fillable in the following lemma.

\begin{lem} \label{cornerlem} If $\Gamma$ is an infinite fillable curve in $\SI$, then $\Gamma$ must be nonoverlapping at the corner.
\end{lem}

\begin{pf} Assume $\Sigma$ is a minimal surface in $\BHH$ with $\PI\Sigma=\Gamma$. As $\Gamma$ is assumed to be tame, $\Gamma^\pm_c=\Gamma\cap (S^1_\infty \times \{\pm\infty\})$ is closed.
	
We claim that $\Gamma^\pm_c$ does not contain any interval in $\CS$. Without loss of generality, assume $I=(p,q)\subset\Gamma^+_c$ where $I$ is an open interval in $S^1_\infty \times \{+\infty\}$. Let  $I'=[p',q']$ be a closed subinterval of $I$. As $\Gamma$ is a Jordan curve, for sufficiently large $c>0$, the region $\widehat{\R}=I'\times[c,\infty)\subset \Si$ is disjoint from $\Gamma$.

For $t>\pi$, let $\R_t= \partial \wh{\R}_t$ where $\wh{\R}_t=I'\times [c,c+t]$. By \cite{Co1}, each $\R_t$ bounds a unique area minimizing surface $\T_t$ with $\PI \T_t=\R_t$ for $t\in(\pi,\infty)$. Note that as $h(\R_t)=t$ is the height of the rectangle $\R_t$, then $h(\R_t)\nearrow\infty$ when $t\nearrow \infty$. 

Now, let $\wh{\R}=[\theta_1,\theta_2]\times[c_1,c_2]$ where $c_2-c_1>\pi$, and let $\R=\partial\wh{\R}$, a finite tall rectangle. Let $\T$ be the unique area minimizing surface in $\BHH$ with $\PI \T=\R$. We claim that if $\Gamma\cap \wh{\R}=\emptyset$, then $\Sigma\cap\T=\emptyset$. Assume not. Let $\Delta$ be the open region in $\BHH$ separated by $\T$, where $\PI\Delta=int(\wh{\R})$. We claim that we can foliate the region $\Delta$ by minimal planes. 

Fix $\theta_3=\dfrac{\theta_1+\theta_2}{2}$. Let $d_1=\dfrac{c_1+c_2-\pi}{2}$ and $d_2=d_1+\pi$. Let $\R'_s=\partial ([\theta_1+(\theta_3-\theta_1)s,\theta_2-(\theta_2-\theta_3)s] \times[c_1+(d_1-c_1)s,c_2-(c_2-d_2)s])$ for $s\in[0,1]$. In particular, $\R'_0=\R$ and $\R'_1=\{\theta_3\}\times[d_1,d_2]$. Notice that, $\{\R'_s \mid s\in(0,1)\}$ foliates $\wh{\R}-\R'_1$. Consider $\{\T'_s\mid s\in(0,1)\}$ where $\T'_s$ is the unique area minimizing surface with $\PI\T'_s=\R'_s$. By \cite{Co1}, $\{\T'_s\mid s\in(0,1)\}$ foliates $\Delta$.

Now, since $\Gamma\cap \wh{\R}=\emptyset$, if $\Sigma\cap\T\neq \emptyset$, $\T$ must separate a compact piece $S$ from $\Sigma$, i.e. $S=\overline{\Delta}\cap\Sigma$. Let $s_0=\sup\{s\in[0,1] \mid S\cap\T'_s\neq \emptyset\}$. By construction, $s_0\in (0,1)$. Then, $\T'_{s_0}$ intersects $S$ tangentially lying in one side of $S$. This contradicts to the maximum principle. This shows that if $\Gamma\cap \wh{\R}=\emptyset$, then $\Sigma\cap\T=\emptyset$.

Hence, for any $t\in(\pi,\infty)$, $\Sigma\cap\T_t=\emptyset$ as $\Gamma\cap\wh{\R}_t=\emptyset$. Recall the infinite rectangles from Section \ref{infcurvesec}. As $h(\R_t)\to \infty$ when $t\to \infty$, by \cite[Prop. 2.1]{ST}, $\T_t$ converges to an infinite rectangle $\widetilde{\T}$. In particular, here $t\to\infty$ corresponds to $d\searrow 1$ case in the proof of \cite[Prop. 2.1]{ST}. Note that $\PI \widetilde{\T}=\widetilde{\R}$ is a union of a pair of vertical line segments $\{p',q'\}\times[c,\infty]$, and a horizontal line segment $[p',q']\times\{c\}$ in $S^1_\infty\times\overline{\BR}$, and a geodesic segment $\gamma'$ in $\BH^2\times\{+\infty\}$ where $\PI\gamma'=\{p',q'\}$.
	
By assumption, $I'$ is in $\Gamma^+_c$. This implies $\PI\Sigma\supset I'$. However, by construction $\T_t\cap \Sigma=\emptyset$ for any $t\in (\pi,\infty)$. Again by construction $\widetilde{\T}\cap\Sigma\neq \emptyset$ as $\widetilde{\R}$ separates $\Gamma$ in $\SI$. Hence, as $\T_t\to \widetilde{\T}$ when $t\to \infty$, for sufficiently large $t_0$, $\T_{t_0}\cap\Sigma\neq \emptyset$. This is a contradiction as we showed that $\Sigma\cap \T_t=\emptyset$ for any $t>\pi$. This shows that $\Gamma$ can not contain any interval in $S^1_\infty \times \{\pm\infty\}$. The proof follows.	
\end{pf}

We can summarize our results in the following corollary:

\begin{cor} \label{fillcor} Let $\Gamma$ be an infinite fillable curve in $\SI$. Then, $\Gamma^\pm_g$ must be a (possibly empty) collection of geodesics, and $\Gamma$ must be nonoverlapping at the corner.
\end{cor}
	
In the next part, we will give some examples of non-fillable infinite curves which holds the conditions in the corollary above.

\subsection{Examples of Non-fillable Curves} \label{nonfilsubsec} \

\vspace{.2cm}

After giving some necessary conditions for infinite fillable curves, the natural question is the following: 

\begin{question}
Are there non-fillable curves $\Gamma$ in $\SI$ which holds the conditions in the Corollary \ref{fillcor}, and with no thin tail?
\end{question}

The answer to this questions is "Yes". There exists non-fillable infinite curves which holds the natural conditions given in Corollary \ref{fillcor}, i.e. being nonoverlapping at the corner, no thin tail, and $\Gamma^\pm_g$ being geodesics.

\begin{figure}[b]
	\begin{center}
		$\begin{array}{c@{\hspace{.4in}}c}

		\relabelbox  {\epsfysize=2in \epsfbox{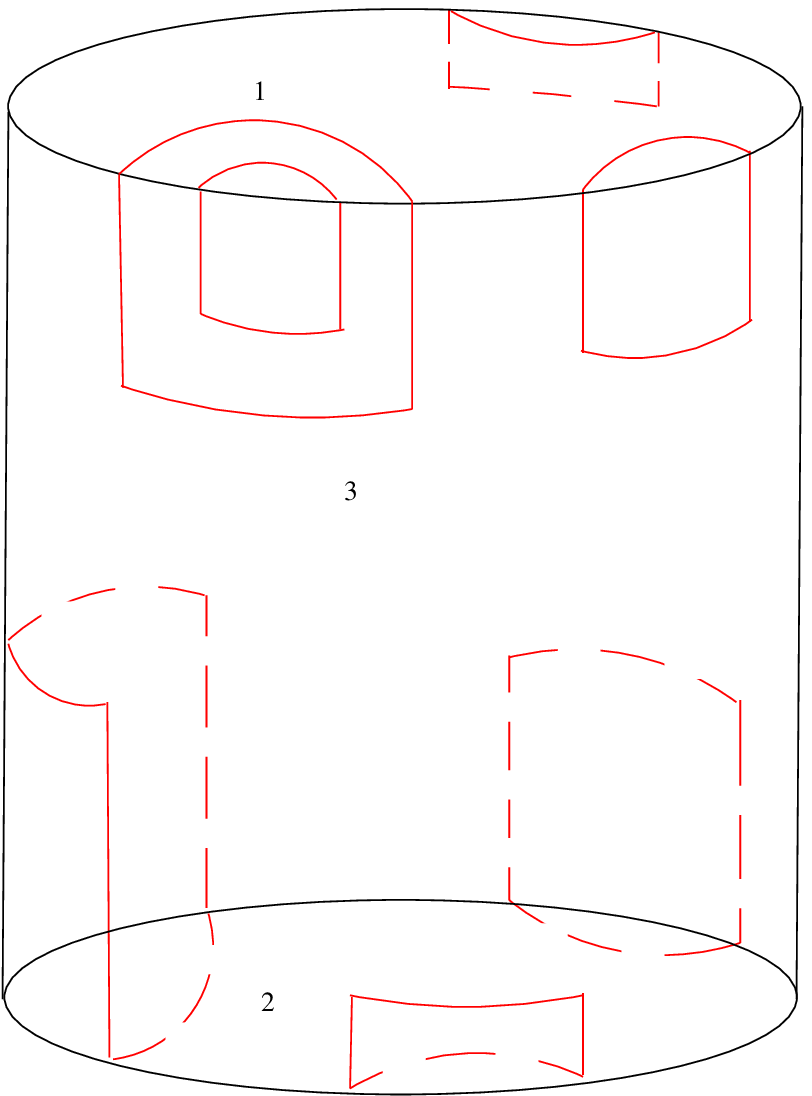}} \relabel{1}{\footnotesize $\Gamma^+$} \relabel{2}{\footnotesize $\Gamma^-$} \relabel{3}{\footnotesize $\R_2^+$} \endrelabelbox &
		
		\relabelbox  {\epsfysize=2in \epsfbox{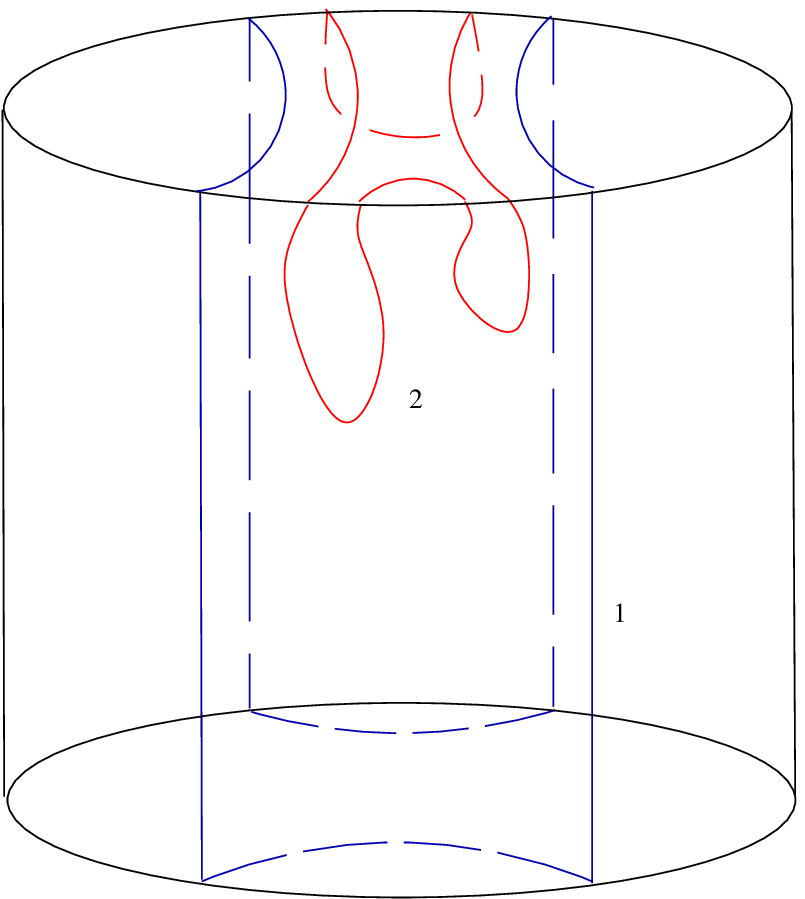}}
		\relabel{1}{\footnotesize $\xi$} \relabel{2}{\footnotesize $\Gamma$ }     \endrelabelbox \\
	\end{array}$
	
\end{center}

\caption{ \label{trappedfig} \footnotesize In the figure left, given $\Gamma^+=\gamma_1^+\cup..\cup\gamma^+_4$ and $\Gamma^-=\gamma^-_1\cup..\cup\gamma^-_3$, we construct fillable $\Gamma$ with infinite rectangles $\R^\pm_i$ for each $\gamma^\pm_i$. In the figure right, $\xi$ is a Scherk curve, and $\Gamma$ is trapped by $\xi$.}
\end{figure}

The examples are called {\em trapped curves}. These are basically the connected infinite curves lying in one side of a Scherk curve $\xi$ (See Figure \ref{trappedfig}). Scherk Curves are very special infinite curves bounding vertical translation invariant Scherk graphs $\{\Sigma_t\}$ \cite{CR}. These are minimal graphs over ideal $2n$-gons in $\BH^2$ where the graph takes values $+\infty$ and $-\infty$ on alternating sides. See \cite{CR, Co2} for further details.

In particular, if $\Sigma$ is a Scherk graph with $\PI\Sigma=\xi$ a Scherk curve, then $\Sigma_t=\Sigma+ t$ is also a Scherk graph with $\PI\Sigma_t=\xi$. Here, $\Sigma+t$ is vertical translation of $\Sigma$ by $t$ for any $t\in\BR$. Let $\Gamma$ be a infinite curve which is one component of $\SI-\xi$ as in the Figure \ref{trappedfig}-right. Let $S$ be a minimal surface in $\BHH$ with $\PI S=\Gamma$. Then, for $s<<0$, $\Sigma_s\cap S=\emptyset$. Hence, for $s_0=\inf \{s\mid \Sigma_s\cap S\neq \emptyset \}$, $\Sigma_{s_0}$ and $S$ will give a contradiction by the maximum principle. See \cite[Theorem 4.6]{Co2} for details of the proof. In particular, we have the following:

\begin{lem} \cite{Co2} There exists  non-fillable infinite curves $\Gamma$ in $\SI$ such that 
		\begin{itemize}
			\item $\Gamma$ is non-overlapping at corner,
			\item $\Gamma^\pm_g$ is a union of geodesics,
			\item $\wt{\Gamma}$ contains no thin tail.
		\end{itemize}	
\end{lem}

\begin{rmk} As it can be seen, trapped curves gives a large family of infinite curves which are non-fillable. Many potential properties to identify infinite fillable curves can be tested by this family. Notice that a connected trapped curve $\Gamma$ is automatically {\em skinny at infinity} by construction (see \cite{Co2}). We will discuss this observation in the next section.
\end{rmk}

\section{Final Remarks} \label{finalsec}

In this part, we will discuss further questions, and remarks. 

\vspace{.2cm}

\subsection{Classification of Fillable Curves} \

Recall that we call $\Gamma$ in $\Si$ \textit{strong fillable} if it bounds an \textit{area minimizing surface}, just \textit{fillable} if it bounds a \textit{minimal surface} in $\BHH$. In \cite{Co1}, we gave the classification of finite strong fillable curves. In \cite{Co2}, we gave the classification of infinite strong fillable curves. However, the examples in this paper shows that fillability question is very different, and the strong fillability results are not enough to answer this question for both infinite and finite curve cases.

The examples of this paper shows that "no thin tail" condition is not enough to classify finite fillable curves in $\Si$. Main difference between strong fillability and fillability is that when $\Gamma$ has more than one component, bounding area minimizing surface for $\Gamma$ still enforce a strong restriction on the position and proximity of all components of $\Gamma$. However, for fillability question, if $\Gamma$ has more than one component, the proximity of the components becomes irrelevant as long as each component is fillable.  That's why fillability question is mainly about fillability of connected curves. Hence, being tall, or similar restriction for non-connected curves becomes useless for fillability question. So, for this problem, it is reasonable to study first the following "simpler" question.

\begin{question} Classify nullhomotopic, connected, fillable curves in $\Si$.
\end{question}


While there were short fillable curves before this paper, our examples give the first non-fillable curves with no thin tail. So, before this paper, we already knew that generalized butterfly curves are not strongly fillable. With the results of this paper, we see that these are also non-fillable curves. On the other hand, the existence of fillable butterfly curves demonstrates that how delicate the classification problem is. So, the following specific class would be a good starting point to understand the classification of finite fillable curves.

\begin{question} Classify all fillable generalized butterfly curves in $\Si$.	
\end{question}

In \cite{KM, Co1}, some examples of fillable butterfly curves given. In this paper, we give examples of non-fillable butterfly curves. After a closer look to both cases, one can notice the following pattern. Let $m$ be the shortest vertical segment in the neck region. Let $M$ be the largest vertical segment in the tall region. Roughly, fillable examples are the ones $m\cdot M> C_0$, and our non-fillable examples are the ones with $m\cdot M< C_1$. However, all these examples are for the simple cases. If the butterfly curve is not horizontal, and the whole curve is very irregular, the current tools are inadequate to finish this delicate problem.

\vspace{.2cm}

\subsection{Con-Fillable Curves:} \

Another interesting version of the asymptotic Plateau problem in $\BHH$ is whether a given collection of curves bounds a connected minimal surface. We will call a curve $\Gamma$ in $\SI$ {\em con-fillable}, if there is a connected filling minimal surface, i.e. $\Gamma$ is con-fillable if there exists a connected, complete, embedded minimal surface $\Sigma$ in $\BHH$ with $\PI \Sigma=\Gamma$.

\vspace{.2cm}

\begin{question} [Con-fillability] \label{conq} Which finite or infinite curves in $\SI$ bound a connected, complete, embedded minimal surface in $\BHH$ ?
\end{question}

In \cite{FMMR}, the authors studied a special case of this problem. They considered a pair of graphical curves $\gamma^\pm$ in $\Si$ where $0<u^+(\theta)-u^-(\theta)<\pi$. They show that if $\gamma^+\cup\gamma^-$ has some rotational symmetry, then it bounds a minimal annulus in $\BHH$. On the other hand, in the same paper, they gave a very neat construction of a non-example for a con-fillable curve for annulus case by using Alexandrov reflection principle \cite[Proposition 5.6]{FMMR}. In particular, they construct a pair of graphical curves $\gamma^+\cup\gamma^-$ which are monotonically getting far from each other, where $\Gamma$ does not bound a minimal annulus.

\subsection{Classification of Infinite Fillable Curves:} \

As discussed in finite curve case, one of the main differences between "strong fillability question" vs. "fillability question" is the number of components. Hence, the key problem in this case is the following:

\begin{question}  \label{conq2} Classify all connected, infinite, fillable curves in $\SI$.
\end{question}

Of course, it is possible a finite collection of curves might be fillable even though each of the curves are not fillable separately. However, considering the difficulty of the problem, above question is a good starting point to study the fillability question for infinite curves.

This classification problem is more complicated than the finite case, as it also contains the same issues in the finite part, and more. One can easily modify an infinite curve to have some butterfly component in the finite part of the curve. So, to separate the difficulties of the finite case, and infinite case, we simplify the question by imposing tall condition.

\begin{question}  \label{conq2} Classify all connected, tall, infinite, fillable curves in $\SI$.
\end{question}

There are various examples of strongly fillable, fillable and non-fillable curves given in \cite{Co2}. One particular non-fillable example is a good starting point to attack the above problem: \textit{Trapped non-fillable curves} (Section \ref{nonfilsubsec}).

Notice that  if $\Gamma$ was fat at infinity, and other expected conditions hold in Lemma \ref{stronglem}, we would know $\Gamma$ would be fillable (indeed strongly fillable). On the other hand, if $\Gamma$ is trapped by a Scherk curve $\xi$, then it is automatically {\em skinny at infinity}. So, trapped curves gives a nice family of examples of non-fillable curves. However, for a given \textit{'skinny at infinity'} $\Gamma^\pm$, by playing with the finite part $\wt{\Gamma}$ of $\Gamma$, it can be arranged so that $\Gamma$ is not trapped by any Scherk curve. So, understanding whether such a skinny at infinity example exist or not would be a key step to classify all such infinite curves in $\SI$. Now, we pose the following conjecture:

\begin{conj} Let $\Gamma$ be an infinite, connected, tall Jordan curve in $\SI$. Then, $\Gamma$ is fillable if and only if $\Gamma$ is fat at infinity.	
\end{conj}

A final remark is about the tame and wild curves in $\SI$. Recall that we call a curve tame if $\Gamma^\pm$ has finitely many component. Otherwise, we will call $\Gamma$ {\em a wild curve.} In this paper, we only considered {\em tame curves} in $\SI$. Interesting examples of fillable, and non-fillable wild curves were constructed in \cite[Section 5.3]{Co2}.

\end{document}